\documentclass{article}

\usepackage{amsfonts}
\usepackage{amssymb}
\usepackage{amsthm}
\usepackage{amsmath}
\usepackage{epsfig}
\usepackage{psfrag}
\usepackage{url}
\usepackage{mathdots}
\usepackage{xcolor}

\def\@begintheorem#1#2{\par\bgroup{\sc #1\ #2. }\it\ignorespaces}
\def\@opargbegintheorem#1#2#3{\par\bgroup{\sc #1\ #2\ (#3). } \it\ignorespaces}
\def\@endtheorem{\egroup}
\newtheorem{theorem}{Theorem}[section]
\newtheorem{corollary}[theorem]{Corollary}
\newtheorem{lemma}[theorem]{Lemma}

\newtheorem{*theorem}[theorem]{*Theorem}
\newtheorem{*corollary}[theorem]{*Corollary}
\newtheorem{*lemma}[theorem]{*Lemma}
\newtheorem{*proposition}[theorem]{*Proposition}

\newcommand{\nofootnote}[1]{}

\def\R{\mathbb{R}}
\def\real{\mathbb{R}}

\def\conv{\operatorname{conv}}
\def\diam{\operatorname{diam}}
\def\ast{\operatorname{ast}}
\def\lk{\operatorname{lk}}

\def\wed{\operatorname{W}}
\def\ops{\operatorname{S}}

\begin{document}

\title{Companion to \\``An update on the Hirsch conjecture''
%: \\Fifty-two years later
}

\author{
Edward D. Kim\thanks{Supported in part by the Centre de Recerca Matem\`atica, NSF grant DMS-0608785 and NSF VIGRE grants DMS-0135345 and DMS-0636297.}
\and
Francisco Santos\thanks{Supported in part by the Spanish Ministry of Science through grant MTM2008-04699-C03-02}
}

\date{}

\maketitle

\begin{abstract}
This is an appendix to our paper ``An update of the Hirsch Conjecture'', containing proofs of some of the results and comments that were omitted in it.
\end{abstract}

\section{Introduction}
This is an appendix to our paper ``An update of the Hirsch Conjecture''~\cite{Kim-Santos-update}, containing proofs of some of the results and comments that were omitted in it. The numbering of sections and results is the same in both papers, although not all appear in this companion. The same occurs with the bibliography, which we repeat here completely although not all of the papers are referenced. The numbering of figures, however,  is correlative. Figures~1 to~6 are in~\cite{Kim-Santos-update} and Figures~7 to~16 are here.

\setcounter{figure}{6}
\setcounter{section}{1}
\section{Bounds and algorithms}
\label{sec:bounds-algorithms}
%\section{Bounds and algorithms}
%\label{sec:positive-results}

\setcounter{subsection}{0}
\subsection{Small dimension or few facets}\label{sec:small-dimension}

\setcounter{theorem}{0}
\begin{*theorem} [Klee~\cite{Klee:PathsII}]
\label{thm:hirschford3}
$H(n,3) = \lfloor \frac{2n}{3} \rfloor - 1$.
\end{*theorem}

\begin{proof} To prove the lower bound, we work in the dual setting where our polytope $P$ simplicial and we want to move from one facet to another along the ridges of $P$.
Figure~\ref{d3lowerbound} shows the graph of a simplicial 3-polytope with nine vertices in which five steps are needed to go from the interior triangle to the most external one (the outer face in the picture, which represents a facet in the polytope). The reader can easily generalize the figure to any number of vertices divisible by  three, adding layers of three vertices that increase the diameter by two. For a number of vertices equal to one or two modulo three, simply add one or two vertices in the interior of the central triangle, subdividing it into three or five triangles. One vertex will not increase the diameter, but two vertices will increase it by one.

\begin{figure}[hbt]
  \begin{center}
    \includegraphics[scale=0.70]{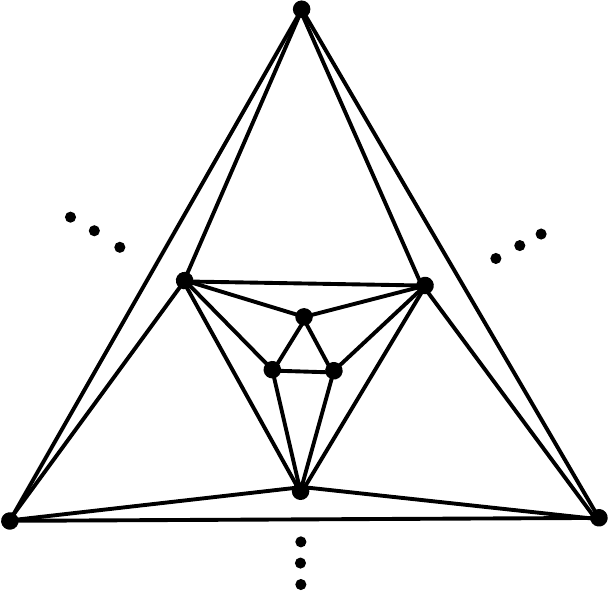}
    \caption{Construction of simplicial $3$-dimensional polytopes with $3+3k$ vertices and diameter $2k+1$. Two steps are needed to cross each layer of skinny triangles.}
       \label{d3lowerbound}
  \end{center}
\end{figure}

For the upper bound, we switch back to simple polytopes.
By double-counting, two times the number of edges of a simple $3$-polytope $P$ equals three times its number of vertices. This, together with Euler's formula, implies that $P$ has exactly $2n-4$ vertices.
Let $u$ and $v$ be two of them. Graphs of 3-polytopes are 3-connected
(see~\cite{Balinski:On-the-graph-structure}, or~\cite{Ziegler:LecturesPolytopes}), which means that 
there are three disjoint paths going from $u$ to $v$. Since the number of intermediate vertices available for these three paths to use is $2n-6$, the shortest of them uses at most $\lfloor \frac{2n}{3}\rfloor-2$ vertices, hence it has at most $\lfloor \frac{2n}{3} \rfloor - 1$ edges.
\end{proof}

It is easy to generalize the second part in the proof to arbitrary dimension, giving the following lower bound. Observe that the formula gives the exact value of $H(n,d)$ for $d=2$ as well.

\setcounter{theorem}{3}
\begin{*proposition}
\[
H(n,d) \ge \left \lfloor\frac{d-1}{d} n \right\rfloor - (d-2).
\]
\end{*proposition}

\begin{proof}
The addition of layers used in the proof of~\ref{thm:hirschford3} can also be described as glueing copies of an octahedron to an already constructed simplicial $3$-polytope. The gluing is along a triangle, so three new vertices are obtained. Before gluing, a projective transformation is made to the octahedron so that the triangle glued is much bigger than the opposite one, which guarantees convexity of the construction.

The generalization to arbitrary dimension is done gluing a  \emph{cross-polytope}, the polar of a $d$-cube. A cross-polytope is also the common convex hull in $\R^d$ of two parallel $(d-1)$-simplices opposite to one another. It has $2d$ vertices and to go from a facet to the opposite one $d$ steps are needed. When the cross-polytope is glued to a given polytope its diameter grows by $d-1$, essentially for the same reasons that will make the proof of Theorem~\ref{thm:monotone} work.
\end{proof}

\setcounter{subsection}{1}
\subsection{General upper bounds on diameters}\label{sec:upper-bounds}

The proof we offer for Theorem~\ref{thm:linear-in-fixed-d} is casically taken from Eisenbrand, H\"ahnle and Rothvoss~\cite{Eisenbrand:limits-of-abstraction}.

\setcounter{theorem}{4}
\begin{*theorem}[Larman~\cite{Larman}]
\label{thm:linear-in-fixed-d}
For every $n>d\ge 3$, $H(n,d)\le n 2^{d-3}$.
\end{*theorem}

\begin{proof}
The proof is by induction on $d$. The base case $d=3$ was Theorem~\ref{thm:hirschford3}.

Let $u$ be an initial vertex of our polytope $P$, of dimension $d>3$. For each other vertex $v\in \operatorname{vert}(P)$ we consider its distance $d(u_1,v)$, and use it to construct a sequence of facets $F_1,\dots,F_k$ of $P$ as follows:

\begin{itemize}
\item Let $F_1$ be a facet that reaches ``farthest from $u$'' among those containing $u$. That is, 
let $\delta_1$ be the maximum distance to $u$ of a vertex  sharing a facet with $u$, and let $F_1$ be that facet. 

\item Let $\delta_2$ be the maximum distance to $u$ of a vertex  sharing a facet with some vertex  at distance $\delta_1+1$ from $u$, and let $F_2$ be that facet. 
%Observe that no assumption is made on the distance from $v_1$ to $u_2$.

\item Similarly, while there are vertices at distance $\delta_i + 1$ from $u$,  let $\delta_{i+1}$ be the maximum distance to $u$ of a vertex  sharing a facet with some vertex  at distance $\delta_i+1$ from $u$, and let $F_{i+1}$ be that facet.
\end{itemize}

We now stratify the vertices of $P$ according to the distances $\delta_1,\delta_2,\dots,\delta_k$ so obtained. 
Observe that $\delta_k$ is the diameter of $P$. By convention, we let $\delta_0=-1$:
\[
V_i:=\{v\in \operatorname{vert}(P) : d(u,v) \in (\delta_{i-1}, \delta_{i}]\}.
\]
We call a facet $F$ of $P$ \emph{active in $V_i$} if it contains a vertex of $V_i$. The crucial property that our stratification has is that no facet of $P$ is active in more than two $V_i$'s. Indeed, each facet is active only in $V_i$'s with consecutive values of $i$, but a facet intersecting $V_i$, $V_{i+1}$ and $V_{i+2}$ would contradict the choice of the facet $F_{i+1}$. In particular, if $n_i$ denotes the number of facets active in $V_i$ we have 
\[
\sum_{i=1}^k n_i \le 2n.
\]

Since each $F_i$ has vertices with distances to $u$ ranging from at least $\delta_{i-1}+1$ to $\delta_i$, 
we have that $\diam(F_i)\ge \delta_{i} - \delta_{i-1} -1$. Even more,
let $Q_i$, $i=1,\dots,k$ be the polyhedron obtained by removing from the facet-definition of $F_i$ the equations of facets of $P$ that are not active in $V_i$ (which may exist since $F_i$ may have vertices in $V_{i-1}$). By an argument similar to the one used for the polyhedron $Q$ of the previous proof, $Q_i$ has still diameter at least $\delta_{i} - \delta_{i-1} -1$. But, by inductive hypothesis, we also have that the diameter of $Q_i$ is at most  $2^{d-4} (n_i-1)$, since it has dimension $d-1$ and at most $n_i-1$ facets.
Putting all this together we get the following bound for the diameter  $\delta_k$ of $P$:

\begin{eqnarray*}
\delta_k &=& \sum_{i=1}^k (\delta_i - \delta_{i-1} -1)  +(k-1) \\
&< & \sum_{i=1}^k 2^{d-4} (n_i - 1) + k \\
&= & 2^{d-4} \sum_i n_i  - k (2^{d-4} -1)
\le 2^{d-3} n.
\end{eqnarray*}

\end{proof}

\setcounter{theorem}{5}
\begin{*theorem}[Kalai-Kleitman~\cite{Kalai:Quasi-polynomial}]
\label{thm:quasipolynomial}
For every $n>d$, $H(n,d)\le n^{\log_2(d)+1}$.
\end{*theorem}

\nofootnote{Here we give a slightly expanded version of the proof:}
\nofootnote{What do we mean with ``slightly expanded''?}

\begin{proof}[Proof of Theorem~\ref{thm:quasipolynomial}]
Let $P$ be a $d$-dimensional polyhedron with $n$ facets, and let
$v$ and $u$ be two vertices of $P$. Let $k_v$ (respectively $k_u$) be the maximal positive number such that the union of all vertices in all paths in $G(P)$ starting from $v$ (respectively $u$) of length at most $k_v$ (respectively $k_u$) are incident to at most $\frac{n}{2}$ facets. Clearly, there is a facet $F$ of $P$ so that we can reach $F$ by a path of length $k_v+1$ from $v$ and a path of length $k_u+1$ from $u$.

We claim that $k_v\leq H_u(\lfloor \frac{n}{2} \rfloor,d)$ (and the same for $k_u$), where $H_u(n,d)$ denotes  the maximum diameter of all $d$-polyhedra with $n$ facets. 
To prove this, let $Q$ be the polyhedron defined by taking only the inequalities of $P$ corresponding to facets that can be reached from $v$ by a path of length at most $k_v$. By construction, all vertices of $P$ at distance at most $k_v$ from $v$ are also vertices in $Q$, and vice-versa. In particular, if $w$ is a vertex of $P$ whose distance from $v$ is $k_v$ then its distance from $v$ in $Q$ is also $k_v$. Since $Q$ has at most $n/2$ facets, we get $k_v\leq H_u(\lfloor \frac{n}{2} \rfloor,d)$.

The claim implies the following recursive formula for $H_u$:
\begin{eqnarray*}
H_u(n,d)
&\le & 2 H_u\left(\left\lfloor \frac{n}{2} \right\rfloor,d\right)  + H_u(n,d-1) + 2,
\end{eqnarray*}
which we can rewrite as
\[
\frac{H_u(n,d)+1}{n}
\le  \frac{H_u\left(\left\lfloor \frac{n}{2} \right\rfloor,d\right)+1}{n/2}  + \frac{H_u(n,d-1)+1}{n}.
\]

This suggests calling $h(k,d):=(H(2^k,d) - 1)/ 2^k$ and applying the recursion with $n=2^k$, to get:
\[
h(k,d)\le h(k-1,d) + h(k,d-1).
\]
This implies $h(k,d)\le \binom{k+d}{d}$, or
\[
H_u(2^k,d) \le 2^k\binom{k + d}{d}.
\]
From this the statement follows if we assume $n\le 2^d$ (that is, $k\le d$). For $n\ge 2^d$ we use Larman's bound 
$H_u(n,d)\le n 2^d \le n^2$, proved below.
\end{proof}

\setcounter{subsection}{3}
\subsection{Some polytopes from combinatorial optimization}\label{sec:special}

\subsubsection*{Small integer coordinates}

\setcounter{theorem}{10}
\begin{*theorem}[Naddef~\cite{Naddef:Hirsch01}] 
\label{thm:01hirsch}
If $P$ is a $0$-$1$ polytope then $\diam(P) \leq n(P) - \dim(P)$.
\end{*theorem}

\begin{proof}
We assume that $P$ is full-dimensional. This is no loss of generality since, if the dimension of $P$ is strictly less than $d$, then $P$ can be isomorphically projected to a face of the cube $[0,1]^d$.

Let $u$ and $v$ be two vertices of $P$. By symmetry, we may assume that $u = (0,\ldots,0)$. If there is an $i$ such that $v_i=0$, then $u$ and $v$ are both on the face of the cube corresponding to $\{{\bf x} \in \R^d \mid x_i = 0\}$, and the statement follows by induction.  Therefore, we assume that $v = (1,\ldots,1)$.
Now, pick any neighboring vertex $v'$ of $v$. There is an $i$ such that $v'_i=0$.  Then, $u$ and $v'$ are vertices of a lower-dimensional $0$-$1$ polytope and we have used one edge to go from $v$ to $v'$. The result follows by induction on $d$.
\end{proof}

\subsubsection*{Transportation and dual transportation polytopes}

We here include the precise definition of $3$-way transportation polytopes, which we skipped in the paper:

\begin{itemize}
\item
{\bf $3$-way axial transportation polytopes.}
Let $a=(a_1,\dots,a_p)$,
$b=(b_1,\dots,b_q)$, and 
$c=(c_1,\dots,c_r)$ be three vectors of lengths $p$, $q$ and $r$, respectively. The $3$-way axial $p \times q \times r$ transportation polytope $P$ given by $a \in \R^p$, $b  \in \R^q$, and $c \in \R^r$ is defined as follows:
\[P = \{(x_{ijk}) \in \R^{p \times q \times r} \mid \sum_{j,k} x_{ijk}=a_i, \sum_{i,k} x_{ijk}=b_j, \sum_{i,j} x_{ijk}=c_k, x_{ijk} \geq 0\}. \]
The polytope $P$ has dimension $pqr-(p+q+r-2)$ and at most $pqr$ facets.

\item
{\bf $3$-way planar transportation polytopes.}
Let $A\in \R^{p\times q}$, $B \in \R^{p \times r}$, and $C \in \R^{q \times r}$ be three matrices.
We define the $3$-way planar $p \times q \times r$ transportation polytope $P$ given by $A$, $B$, and $C$ as follows:
\[P = \{(x_{ijk}) \in \R^{p \times q \times r} \mid \sum_{k} x_{ijk}=A_{ij}, \sum_{j} x_{ijk}=B_{ik}, \sum_{i} x_{ijk}=C_{jk}, x_{ijk} \geq 0\}. \]
It has dimension $(p-1)(q-1)(r-1)$ and at most $pqr$ facets.

\end{itemize}

\setcounter{subsection}{4}
\subsection{A continuous Hirsch conjecture}\label{sec:continuous}

Let us expand a bit the concept of curvature of the central path and its relation to the simplex method.
For further description of the method we refer  to the books~\cite{Boyd:ConvexOptimization, Renegar:A-Mathematical-View}.

The central path method is one of the interior point methods  for solving a linear program.
As in the simplex method, the idea is to move from a feasible point to another feasible point on which the given objective linear functional is improved. In contrast to the simplex method, where the path travels from vertex to neighboring vertex along the graph of the feasibility polyhedron $P$, this method follows a certain curve through the strict  interior of the polytope.

More precisely, to each linear program,
\[
\text{ Minimize } c\cdot {\bf x}, \text{ subject to } A{\bf x} = {\bf b} \text{ and } {\bf x} \ge 0,
\]
the method associates a \emph{(primal) central path}  $\gamma_c : [0, \beta) \rightarrow \R^d$ which is an analytic curve through the interior of the feasible region and such that $\gamma_c(0)$ is an optimal solution of the problem. The central path is well-defined and unique even if the program has more than one optimal solution, but its definition is implicit, so that there is no direct way of computing $\gamma_c(0)$. To get to $\gamma_c(0)$, one starts at any feasible solution and tries to follow a curve that approaches more and more the central path, using for it certain barrier functions. (Barrier functions play a role similar to the choice of pivot rule in the simplex method.  
The standard barrier function is the logarithmic function $f(x) = -\sum_{i=1}^n \ln(A_i x - b_i)$.)  

Of course, it is not possible to follow the curve exactly. Rather, one does Newton-like steps trying not to get too far.
How much can one improve in a single step is related to the curvature of the central path: if the path is rather straight one can do long steps without deviating too far from it, if not one needs to use shorter steps. Thus, 
the  \emph{total curvature} $\lambda_c(P)$ of the central path, defined in the usual differential-geometric way, can be considered a continuous analogue of the diameter of the polytope $P$, or at least of the maximum distance from any vertex to a vertex maximizing the functional $c$.

\setcounter{section}{2}
\section{Constructions}
\label{sec:constructions}

\setcounter{subsection}{0}
\subsection{The wedge operation}\label{sec:wedging}
\begin{figure}[hbt]
\begin{center}
\includegraphics[scale=0.6]{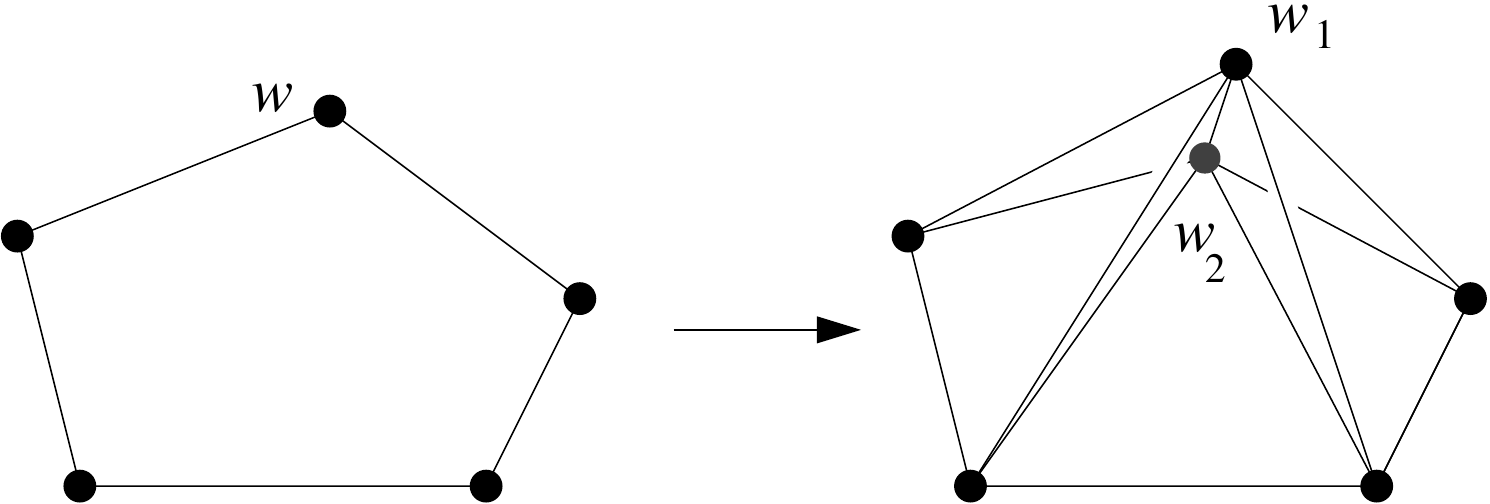}
\caption{The simplicial version of Figure 3 in~\cite{Kim-Santos-update}: a $5$-gon and a one-point suspension on its topmost vertex}\label{fig:ops}
\end{center}
\end{figure}

The dual operation to wedging, usually performed for simplicial polytopes (or for simplicial complexes in general), is the \emph{one-point suspension}. We refer the reader to~\cite[Section 4.2]{triang-book} for an expanded overview of this topic. Let $w$ be a vertex of the polytope $P$. The one-point suspension of $P\subset \real^d$ at the vertex $w$ is the polytope
\[
\ops_w(P) := \conv\big( (P \times \{0\}) \cup (\{w\} \times \{-1,+1\})\big)\subset\real^{d+1}.
\]
That is, $\ops_w(P)$ is formed by taking the convex hull of $P$ (in an ambient space of one higher dimension) with a ``raised'' and ``lowered'' copy of the vertex $w$. See Figure~\ref{fig:ops} for an example.

Recasting Lemma~{3.1} to the dual setting gives the following simplicial version of it:

\setcounter{theorem}{0}
\begin{lemma}\label{lemma:ops}
Let $P$ be a $d$-polytope with $n$ vertices. Let $P'=\ops_w(P)$ be its one-point suspension on a certain vertex $w$. Then $P'$ is a $(d+1)$-dimensional polytope with $n+1$ vertices, and the diameter of the dual graph of $P'$ is at least the diameter of the dual graph of $P$.
\end{lemma}

The one-point suspension of a simplicial polytope is a simplicial polytope. In fact, 
the one-point suspension can be described at the leval of 
 abstract simplicial complexes:
Let $L$ be a simplicial complex and $w$ a vertex of it. Recall that the anti-star $\ast_L(w)$ of $w$ is the subcomplex consisting of simplices not using $w$ and the link $\lk_L(w)$ of $w$ is the subcomplex of simplices not using $w$ but joined to $w$. If $L$ is a PL $k$-sphere, then $\ast_L(w)$  and $\lk_L(w)$ are a $k$-ball and a $(k-1)$-sphere, respectively. The one-point suspension of $L$ at $w$ is the following complex:
\[
\ops_w(L):=(\ast_L(w) * w_1)\cup (\ast_L(w) * w_2) \cup (\lk_L(w)* \overline{w_1w_2}).
\]
Here $*$ denotes the \emph{join} operation: $L*K$ has as simplices all joins of one simplex of $K$ and one of $L$. In Figure~\ref{fig:ops} the three parts of the formula are the three triangles using $w_1$ but not $w_2$, the three using $w_2$ but not $w_1$, and the two using both, respectively.

In Section~\ref{sec:hirsch-sharp-polytopes}
we will make use of an iterated one-point suspension. That is, in $\ops_w(P)$ we take the one-point suspension over one of the new vertices $w_1$ and $w_2$, then again in one of the new vertices created, and so on. We leave it to the reader to check that, at the level of simplicial complexes, the one-point suspension iterated $k$ times produces the following simplicial complex, where $\Delta_k$ is a $k$-simplex with vertices $w_1,\dots,w_{k+1}$ and $\partial \Delta_k$ is its boundary. Observe that this generalizes the formula for $\ops_w(L)$ above:
\[
\ops_w(L)^{(k)}:= (\ast_L(w) * \partial \Delta_k) \cup (\lk_L(w)*{\Delta_k}).
\]

\setcounter{subsection}{1}
\subsection{The $d$-step and non-revisiting conjectures}\label{sec:equivalences}

In this section we had proof that for both the Hirsch and the non-revisiting conjectures the general case is equivalent to the case $n=2d$, but we did not finish proving that the two were equivalent:

\setcounter{theorem}{6}
\begin{*theorem}[Klee-Walkup~\cite{Klee:d-step}]
\label{thm:dstep-nonrevisiting}
The Hirsch, non-revisiting, and $d$-step Conjectures~{1.1},~{3.3}, and~{3.6} are equivalent.
\end{*theorem}

\begin{proof}
\nofootnote{I corrected this proof (a second induction was necessary), but am tempted to move it to the appendix. Paco}
Clearly, the $d$-step conjecture is a special case of both the Hirsch and the non-revisiting conjectures. By Theorems~{3.2} and~{3.4}, to prove that the $d$-step conjecture implies the other two we may restrict our attention to polytopes of dimension $d$ and with $2d$ facets. We also use induction on the \emph{codimension}. That is, we assume the Hirsch and non-revisiting conjectures for all polytopes with number of facets minus dimension smaller than $d$.

Let $u$ and $v$ be two vertices of a $d$-polytope $P$ with $2d$ facets. We will also induct on the number of common facets containing both $u$ and $v$. The base case is when $u$ and $v$ are complementary, in which the $d$-step conjecture applied to them gives a non-revisiting path of length at most $d$. 

So, we assume that  $u$ and $v$ are in a common facet $F$ of $P$. $F$ has at most $2d-1$ facets itself.

\begin{itemize}
\item If  $F$ has less than $2d-1$ facets, then $F$ has the non-revisiting and Hirsch properties by induction on ``number of facets minus dimension'', and we are done.

\item If $F$ has $2d-1$ facets, since it has dimension $d-1$ there is a facet $G$ of $F$ not containing $u$ nor $v$. Let $P'=\wed_{G}(F)$ be the wedge of $F$ on $G$. Let $u_1$ and $v_2$ be vertices of $P'$ projecting to vertices $u$ and $v$ of $P$ and such that $F_1$ contains $u_1$ and $F_2$ contains $v_2$. As in the proof of Theorem~{3.4}, $F_1$ and $F_2$ denote the non-vertical facets of the wedge $P'$. $P'$ again has dimension $d$ and $2d$ facets, but its vertices $u_1$ and $v_2$ have one less facet in common than $u$ and $v$ had. By induction on the number of common facets, there is a non-revisiting path of length at most $d$ between $u_1$ and $v_2$ in $P'$. When this path is projected to $F$, it retains the non-revisiting property and its length does not increase.
\end{itemize}
\end{proof}

\setcounter{subsection}{2}
\subsection{The Klee-Walkup polytope $Q_4$} \label{sec:Q4}

Let us give further details on the structure of the  Hisrsch-sharp polytope $Q_4$ constructed by Klee and Walkup.
Recall that the coordinates we use for the nine vertices of $Q_4$ are:

\[
\begin{tabular}{lll}
&$w:=(0,0,0,-2)$,\\
$a:= (-3,3,1,2)$,&&$e:= (3,3,-1,2)$,\\
$b:=(3,-3,1,2),$&&$f:=(-3,-3,-1,2),$\\
$c:=(2,-1,1,3)$,&&$g:=(-1,-2,-1,3)$,\\
$d:=(-2,1,1,3)$,&&$h:=(1,2,-1,3)$.\\
\end{tabular}
\]

What follows is the input and output of the polymake~\cite{polymake} computation of the face complex of $Q_4$. The input vertices are given in homogenized version, which means and additional coordinate of 1's is added to each.

\begin{verbatim}
POINTS 
 1  0  0  0 -2
 1 -3  3  1  2  
 1  3 -3  1  2
 1  2 -1  1  3
 1 -2  1  1  3
 1  3  3 -1  2 
 1 -3 -3 -1  2
 1 -1 -2 -1  3
 1  1  2 -1  3
\end{verbatim}

The output {\tt VERTICES\_IN\_FACETS} lists the facets as sets of vertices. Polymake numbers the vertices starting with 0, so our vertices $w,a,\dots,h$ become labeled {\tt 0, 1,\dots,8}:

\begin{verbatim}
VERTICES_IN_FACETS
 {2 3 7 8}
 {0 1 2 3}
 {1 2 3 4}
 {2 3 6 7}
 {2 3 4 6}
 {0 2 4 6}
 {0 2 6 7}
 {0 1 2 4}
 {1 6 7 8}
 {0 1 6 8}
 {1 4 7 8}
 {0 1 4 6}
 {1 4 6 7}
 {3 4 6 7}
 {3 4 7 8}
 {0 5 6 8}
 {5 6 7 8}
 {0 1 5 8}
 {1 4 5 8}
 {3 4 5 8}
 {0 1 3 5}
 {1 3 4 5}
 {0 5 6 7}
 {0 2 5 7}
 {2 5 7 8}
 {0 2 3 5}
 {2 3 5 8}
\end{verbatim}

\begin{figure}[hbt]
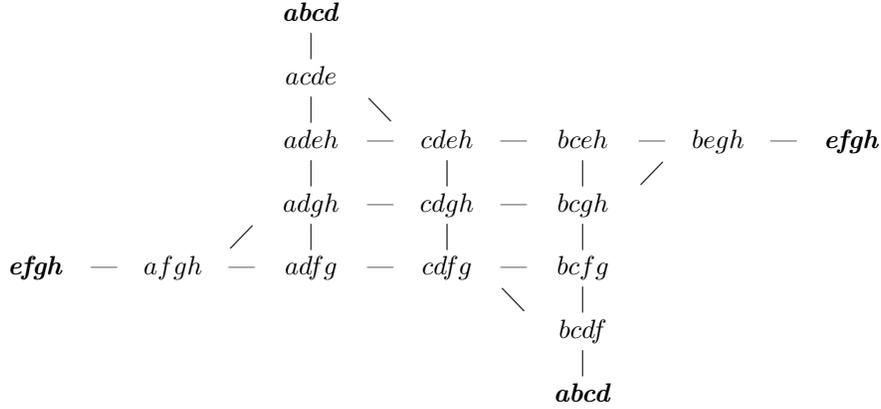

\[
\begin{matrix}
 & & \textbf{\textit{abcd}} && && &&&\\
 && | & & && &&  &\\
 & & acde && && &&&\\
  && | & \diagdown & &&&& &\\
  && adeh &\text{---} & cdeh &\text{---} & bceh &\text{---} & begh&\text{---}\quad  \textbf{\textit{efgh}} \\
  && | && | && | & \diagup & &\\
  && adgh &\text{---} & cdgh &\text{---} & bcgh  && &\\
   & \diagup & | && | && | &&&\\
 \textbf{\textit{efgh}}\quad\text{---}\quad afgh & \text{---} & adfg &\text{---} & cdfg &\text{---} & bcfg &&&\\
   && && & \diagdown & |  &&&\\
  && && & & bcdf & &&\\
   && && & & |  &&&\\
     && && & & \textbf{\textit{abcd}} &  &&\\
 \end{matrix}
\]
\caption{The dual graph of the subcomplex $K$}\label{fig:dualK}
\end{figure}

You should verify that there are exactly 15 tetrahedra not using $w$ (the label {\tt 0}) are precisely the ones in Figure~\ref{fig:dualK}.

\begin{figure}[htb]
\begin{center}
\includegraphics[width=1.5in]{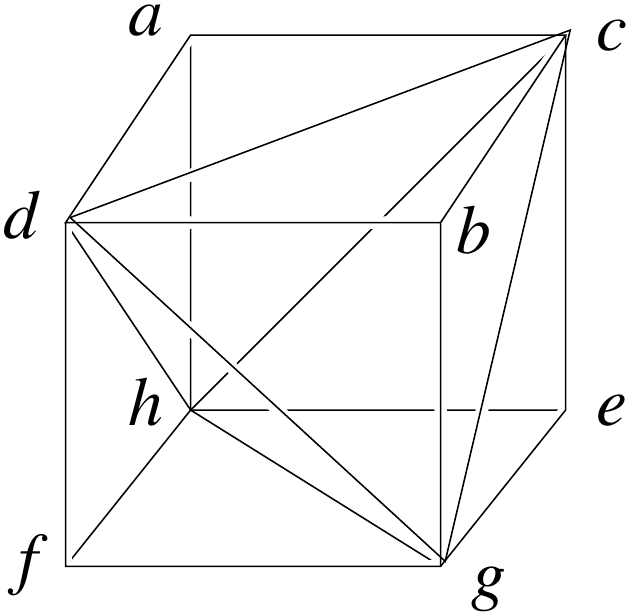}
\caption{The link of $w$ in $Q_4$  is combinatorially a triangulation of the boundary of a cube.}
\label{fig:klee-walkup-cube}
\end{center}
\end{figure}

\begin{figure}[htb]
\begin{center}
\includegraphics[width=4in]{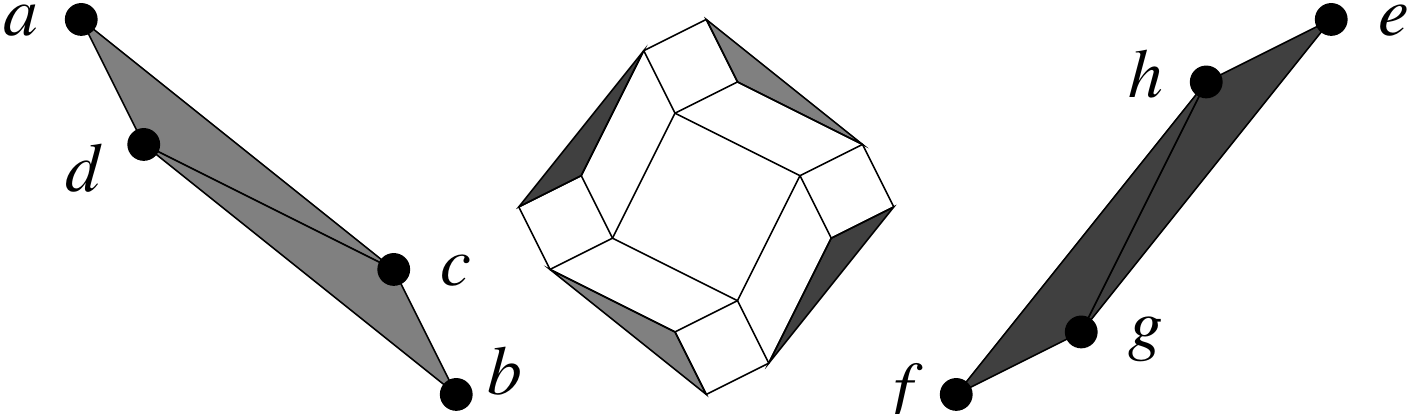}
\caption{The Klee-Walkup complex as a mixed subdivision. The shadowed triangles represent tetrahedra adjacent to $abcd$ and $efgh$}
\label{fig:klee-walkup}
\end{center}
\end{figure}

From the picture we can also read the tetrahedra of $\partial Q_4^*$ that use $w$: there is one for each triangle that appears only once in the list. For example, since $abcd$ is adjacent only to $acde$ and $abcd$,  the triangles $abc$ and $bcd$ are joined to $w$ . The boundary of the antistar of $w$, that is, the \emph{link} of $w$ in $Q_4^*$ turns out to be, combinatorially, the triangulation of the boundary of a cube displayed in Figure~\ref{fig:klee-walkup-cube}.
The anti-star $K$ of $w$ in $\partial Q_4^*$ is a topological triangulation of the interior of the cube. But we need to deform the cube a bit to realize this triangulation geometrically. This is shown in Figure~\ref{fig:klee-walkup}: the quadrilaterals $abcd$ and $efgh$ are displayed separately as lying in two different horizontal planes (so that the two relevant tetrahedra $abcd$ and $efgh$ degenerate to flat quadrilaterals), and the central part of the figure shows the intersection of $K$ with their bisecting plane. Tetrahedra with three points on one plane and one in the other appear as triangles and tetrahedra with two points on either side appear as quadrilaterals. The tetrahedra $abcd$ and $efgh$ do not show up in the figure, since they do not intersect the intermediate plane. For the interested reader, this picture is an example of a \emph{mixed subdivision} of the Minkowski sum of two polygons. The fact that triangulations of polytopes with their vertices lying in two parallel hyperplanes can be pictured as mixed subdivisions is the \emph{polyhedral Cayley trick}~\cite[Chapter 9]{triang-book}.

\setcounter{section}{3}
\subsection{Many Hirsch-sharp polytopes?}\label{sec:hirsch-sharp-polytopes}

\subsubsection*{Trivial Hirsch-sharp polytopes}\label{sec:examples}

\setcounter{theorem}{9}
\begin{*proposition}
For every $n \ge d$ there are simple unbounded $d$-polyhedra with $n$ facets and diameter $n-d$.
\end{*proposition}

\begin{figure}[htb]
\begin{center}
\includegraphics[width=1.7in]{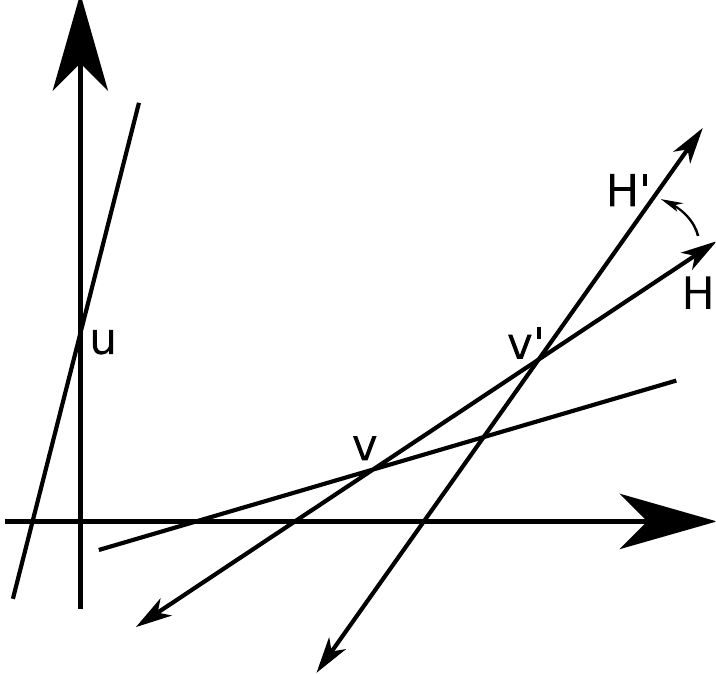}
\caption{Tilting the hyperplane $H$, example in dimension two}
\label{fig:unbounded-hirsch-sharp}
\end{center}
\end{figure}

\begin{proof}
The proof is by induction on $n$, the base case $n=d$ being the orthant $\{x_i\ge 0, \forall i\}$.
Our inductive hypothesis is not only that  we have constructed a $d$-polyhedron $P$ with $n-1$ facets and diameter $n-d-1$; also, that vertices $u$ and $v$ at distance $n-d-1$ exist in it with $v$ incident to some unbounded ray $l$. Let $H$ be a supporting hyperplane of $l$, and tilt it slightly at a point $v'$ in the interior of $l$ to obtain a new hyperplane $H'$. See Figure~\ref{fig:unbounded-hirsch-sharp}. Then, the polyhedron $P'$ obtained cutting $P$ with the tilted hyperplane $H'$ has $n$ facets and diameter $n-d$; $v$ is the only vertex adjacent to $v'$ in the graph, so we need at least $1 + (n-d-1)$ steps to go from $v'$ to $u$.
\end{proof}

\subsubsection*{Non-trivial Hirsch-sharp polytopes}\label{sec:many-hirsch-sharp}

In~\cite{Kim-Santos-update} we only proved part (1) of the following result:

\setcounter{theorem}{10}
\begin{*theorem}[Fritzsche-Holt-Klee~\cite{Fritzsche99morepolytopes,Holt:Hsharpd7,Holt:Many-polytopes}]
\label{thm:hirsch-sharp}
Hirsch-sharp $d$-polytopes with $n$ facets exist in at least the following cases:
%(1) $n\le 2d$;
(1) $n\le 3d-3$; and
(2) $d\ge 7$.
\end{*theorem}

\begin{figure}[hbt]
\begin{center}
\includegraphics[scale=0.60]{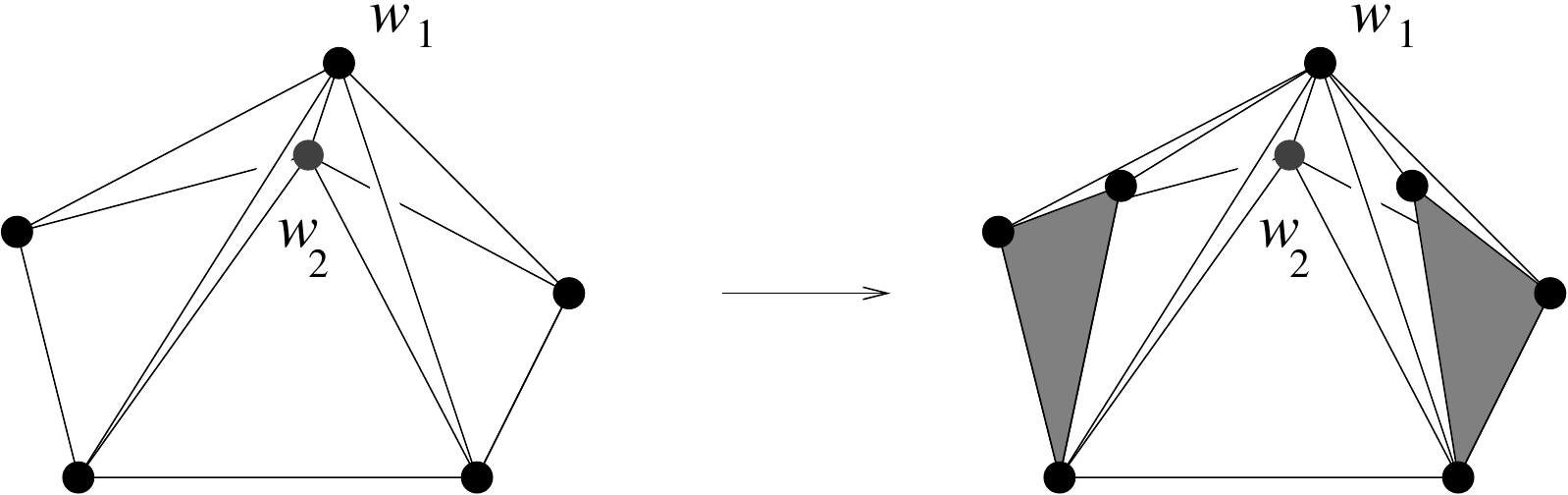}
\caption{The simplicial version of~\cite[Figure~5]{Kim-Santos-update}. Wedging becomes one-point suspension and truncation is stellar subdivision}\label{fig:sharp-3d-3-ops}
\end{center}
\end{figure}

The proof of part (2) is easier to understand in the simplicial framework. So, as a warm-up, we include (see Figure \ref{fig:sharp-3d-3-ops}) the simplicial version of~\cite[Figure~5]{Kim-Santos-update}. We already know that the polar of wedging is one-point suspension. The polar of truncation of a vertex is the \emph{stellar subdivision} of a facet by adding to our polytope a new vertex very close to that facet.

The key property in the proof of Lemma~{3.12} is that the wedge and one-point suspension operations  do not only preserve Hirsch-sharpness; they also increase the number of vertices or facets (respectively) that are at Hirsch distance from one another. This suggests looking at what happens when we iterate the process. The answer, that we state in the simplicial version, is as follows:

\setcounter{theorem}{13}
\begin{lemma}[Holt-Klee~\cite{Holt:Many-polytopes}]
\label{lemma:ops-iterated}
Let $P$ be a simplicial $d$-polytope with more than $2d$ vertices. Let $A$ and $B$ be two facets of it at Hirsch distance in the dual graph and let $w$ be a vertex contained in neither $A$ nor $B$. Let $P^{(k)}$ be the $k^\textrm{th}$ one-point suspension of $P$ on the vertex $w$. 

Then, $P^{(k)}$ has two $(k+1)$-tuples of facets $\{A_1,\dots,A_{k+1}\}$ and $\{B_1,\dots,B_{k+1}\}$ with every $A_i$ at Hirsch distance from every $B_i$. All the facets in each tuple are adjacent to one another.
\end{lemma}

\begin{proof}
We use the following formula, from Section~\ref{sec:wedging}, for the iterated one-point suspension of the simplicial complex $L=\partial P$:
\[
\ops_w(L)^{(k)}:= (\ast_{L}(w) * \partial \Delta_k) \cup (\lk_{L}(w)* {\Delta_k}).
\]
Here $\Delta_k$ is a $k$-simplex. The two groups of facets in the statement are $A * \partial \Delta_k$ and  $B * \partial \Delta_k$. The details are left to the interested reader.
\end{proof}

%This is the basis for the following result of Fritzsche and Holt. 

%\begin{corollary}[Fritzsche-Holt~\cite{Fritzsche99morepolytopes,Holt:Hsharpd7}]
%\label{coro:glue}
%There are Hirsch-sharp $7$-polytopes with any number of facets.
%\end{corollary}

\begin{proof}[Proof of part (2) of Theorem~\ref{thm:hirsch-sharp}.]
We include only the proof for the case $d\ge 8$, contained in~\cite{Fritzsche99morepolytopes}. The improvement to $d=7$ was later found by Holt~\cite{Holt:Hsharpd7}.\nofootnote{Should/could we include a proof of $d=7$?}

Both are based on a new operation on polytopes that we now introduce. The version for simple polytopes is called \emph{blending}, but we describe it for simplicial polytopes and call it \emph{glueing}. Glueing 
is simply a combinatorial/geometric version of the \emph{connected sum} of topological manifolds.
Let $P_1$ and $P_2$ be two simplicial $d$-polytopes and let $F_1$ and $F_2$ be respective facets. The manifolds are $\partial P_1$ and $\partial P_2$ (two $(d-1)$-spheres); from them we remove  the interiors of $F_1$ and $F_2$ after which we glue their boundaries. See Figure~\ref{fig:glueing}, where the operation is performed on two facets of the same polytope. On the top part we glue the polytopes ``as they come'', which does not preserve convexity. But if projective transformations are made on $P_1$ and $P_2$ that send points that are close to  $F_1$ and $F_2$ to infinity, then the glueing  preserves convexity, so it yields a polytope that we denote $P_1 \# P_2$. This is shown on the bottom part of the Figure. 

\begin{figure}[hbt]
\begin{center}
\includegraphics[scale=0.60]{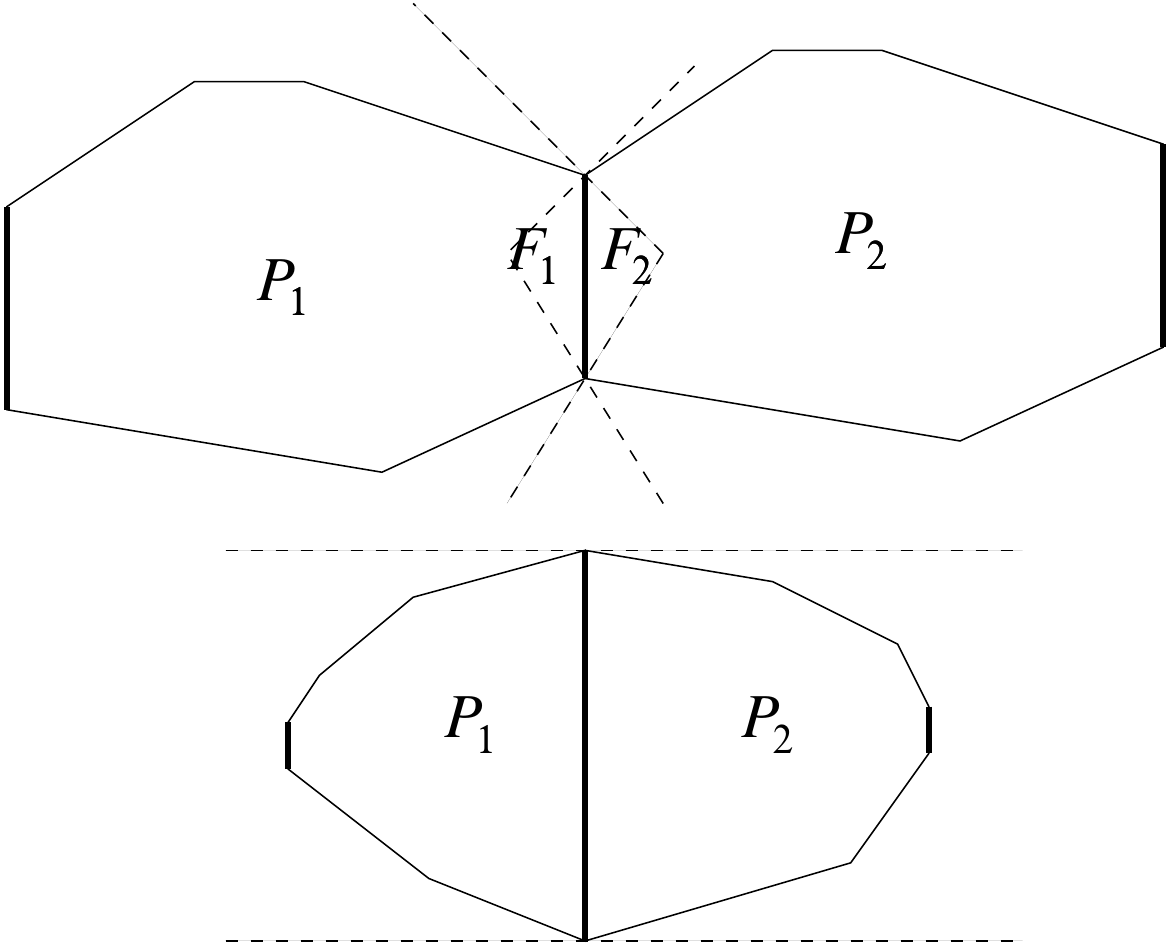}
\caption{Glueing two simplicial polytopes along one facet. In the version on the bottom, a projective transformation is done to $P_1$ and $P_2$ before glueing, to guarantee convexity of the outcome}\label{fig:glueing}
\end{center}
\end{figure}

Glueing \emph{almost} adds the diameters of the two original polytopes. Suppose that the facets $F_1$ and $F_2$ are at distances $\delta_1$ and $\delta_2$ to certain facets $F'_1$ and $F'_2$ of $P_1$ and $P_2$. Then, to go from $F'_1$ to $F'_2$ in $P_1 \# P_2$ we need at least $(\delta_1 -1) + 1 + (\delta_2-1) = \delta_1 + \delta_2 - 1$ steps.

But we can do better if we combine glueing with the iterated one-point suspension.
Consider the simplicial Klee-Walkup $4$-polytope $Q_4^*$ described in Section~\ref{sec:Q4} and let $A$ and $B$ two facets of it at distance five. Let $P'$ be the $4^\textrm{th}$ one-point suspension of it on the vertex $w$ not contained in $A\cup B$. Observe that $P'$ has 13 vertices and dimension eight. 
By the lemma, $P'$ has two groups of five facets $\{A_1,\dots,A_{5}\}$ and $\{B_1,\dots,B_{5}\}$ with every $A_i$ at Hirsch distance from every $B_i$ and all the facets in each group adjacent to one another. 

We now glue several copies of $P'$ to one another, a $B_i$ from each copy glued to an $A_i$ of the next one.  Each glueing adds five vertices and, in principle, four to the diameter. But Lemma~\ref{lemma:ops-iterated} implies the following nice property for $P'$: half of the eight facets adjacent to each $A_i$ are at distance four to half of the facets adjacent to each $B_i$. Using the language of Fritzsche, Holt and Klee, we call those facets the \emph{slow neighbors} of each $A_i$ or $B_i$, and call the others \emph{fast}. Since half of the total neighbors are slow, we can make all glueings so that every fast neighbor is glued to a slow one and vice-versa. This increases the diameter by one at every glueing, and the result is Hirsch-sharp.

The above construction yields Hirsch-sharp 8-polytopes with $13 + 5k$ vertices, for every $k\ge 0$. We can get the intermediate values of $n$ too, via truncation. By Lemma~{3.12}, every time we do a one-point suspension on a Hirsch-sharp simplicial polytope we can increase the number of facets by one or two via a stellar subdivision at each end. Since the polytope $P'$ we are glueing is a 4-fold one-point suspension, and since there are two ends that remain unglued (the $A$-face of the first copy and the $B$-face of the last) we can do up to eight stellar subdivisions to it and still preserve Hirsch-sharpness.
\end{proof}

%The result has been improved to dimension $7$ by Holt, which finishes the proof of part (2) of Theorem~\ref{thm:hirsch-sharp}:

%\begin{corollary}[Holt~\cite{Holt:Hsharpd7}]
%\label{coro:glue2}
%There are Hirsch-sharp $8$-polytopes with any number of facets.
%\end{corollary}

%\begin{proof}
%TO BE DONE
%\end{proof}

%%%%%%%%%%%%%%

\setcounter{subsection}{4}
\subsection{The unbounded and monotone Hirsch conjectures are false}\label{sec:unbounded-monotone}

\setcounter{theorem}{15}
\begin{*theorem}[Todd~\cite{Todd:MonotonicBoundedHirsch}]
\label{thm:monotone}
There is a simple bounded polytope $P$, two vertices $u$ and $v$ of it, and a linear functional $\phi$ such that:
\begin{enumerate}
\item $v$ is the only maximal vertex for $\phi$.
\item Any edge-path from $u$ to $v$ and monotone with respect to $\phi$ has length at least five.
\end{enumerate}
\end{*theorem}

\begin{proof}
Let $Q_4$ be the Klee-Walkup polytope. 
Let $F$ be the same ``ninth facet'' as in the previous proof, one that is not incident to the two vertices $u$ and $v$ that are at distance five from each other. Let $H_2$ be the supporting hyperplane containing $F$ and let $H_1$ be any supporting hyperplane at the vertex $v$. Finally, let $H_0$ be a hyperplane containing the (codimension two)
intersection of $H_1$ and $H_2$ and which lies ``slightly beyond $H_1$'', as in Figure~\ref{fig:monotone}. (Of course, if $H_1$ and $H_2$ happen to be parallel, then $H_0$ is taken to be parallel to them and close to $H_1$.)
 The exact condition we need on $H_0$ is that it does not intersect $Q_4$ and the small, wedge-shaped region between $H_0$ and $H_1$ does not contain the intersection of any 4-tuple of facet-defining hyperplanes of $Q_4$. 

\begin{figure}[htb]
\begin{center}
\includegraphics[width=\textwidth]{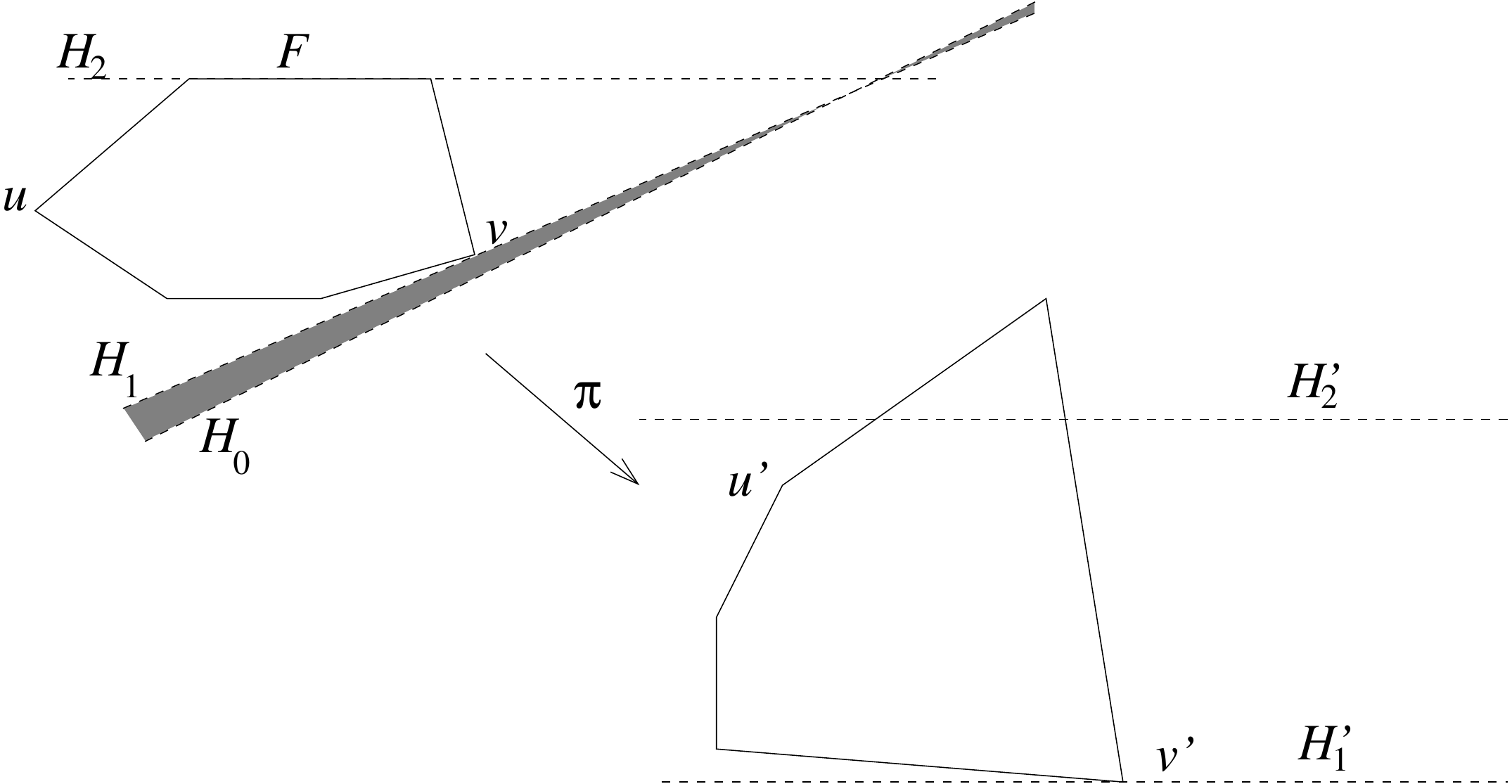}
\caption{Disproving the monotone Hirsch conjecture}
\label{fig:monotone}
\end{center}
\end{figure}

We now make a projective transformation $\pi$ that sends $H_0$ to be the hyperplane at infinity. In the polytope $Q'_4=\pi(Q_4)$ we ``remove'' the facet $F'=\pi(F)$ that is not incident to the two vertices $u'=\pi(u)$ and $v'=\pi(v)$. That is, we consider the polytope $Q''_4$ obtained from
$Q'_4$ by forgetting the inequality that creates the facet $F'$ (see Figure~\ref{fig:monotone} again).  Then $Q''_4$ will have new vertices not present in $Q'_4$, but it also has the following properties:

\begin{enumerate}
\item $Q''_4$ is bounded. Here we are using the fact that the wedge between $H_0$ and $H_1$ contains no intersection of facet-defining hyperplanes: this implies that no facet of $Q''_4$ can go ``past infinity''.

\item It has eight facets: four incident to $u'$ and four incident to $v'$.

\item The functional $\phi$ that is maximized at $v'$ and constant on its supporting hyperplane $H'_1=\pi(H_1)$ is also constant on $H'_2=\pi(H_2)$, and $u'$ lies on the same side of $H'_1$ as $v'$.
\end{enumerate}

In particular, no $\phi$-monotone path from $u'$ to $v'$ crosses $H'_1$,
which means it is also a path from $u'$ to $v'$ in the polytope $Q'_4$, combinatorially isomorphic to  $Q_4$.
\end{proof}

\nofootnote{Should/could we include the proof of Theorem~{3.17}?

In both the constructions of Theorems~{3.13} and~{3.15} one can glue several copies of the initial block $Q_4$ to one another,
%The basic idea is (the dual of) the same one used in Corollary~\ref{coro:glue}. We skip details,
%but in both cases we 
increasing the number of facets by four and the diameter by five, per $Q_4$ glued:

\setcounter{theorem}{16}
\begin{*theorem}[Klee-Walkup, Todd]
\begin{enumerate}
\item There are unbounded $4$-polyhedra with $4+4k$ facets and diameter $5k$, for every $k\ge 1$.

\item There are bounded $4$-polyhedra with $5+4k$ facets and vertices $u$ and $v$ of them with the property that any monotone path from $u$ to $v$ with respect to a certain linear functional $\phi$ maximized at $v$ has length at least $5k$.
\end{enumerate}
\end{*theorem}
%\begin{proof}
%TO BE DONE
%\end{proof}
}

\setcounter{subsection}{5}
\subsection{The topological Hirsch conjecture is false}\label{sec:topological}

\setcounter{theorem}{17}
\begin{*theorem}[Mani-Walkup~\cite{Mani:A-3-sphere-counterexample}]
\label{thm:mani-walkup}
There is a triangulated 3-sphere with 16 vertices and without the non-revisiting property. 
Wedging on it eight times yields a non-Hirsch $11$-sphere with $24$ vertices.
\end{*theorem}

\begin{proof}
The key part of the construction is the two-dimensional simplicial complex $K$ consisting of the following 32 triangles:
\[
\begin{tabular}{cccccccc}
amr&mbr&bnr&ncr&cor&odr&dpr&par\\
amt&mbt&bnt&nct&cot&odt&dpt&pat\\
\end{tabular}
\]
%and
\[
\begin{tabular}{cccccccc}
aoq&obq&bpq&pcq&cmq&mdq&dnq&naq\\\
aos&obs&bps&pcs&cms&mds&dns&nas\\
\end{tabular}
\]
The first and second halves are topological 2-spheres, triangulated in the form of 
double pyramids over the octagons $ambncodp$ and $aobpcmdn$ (same vertices, but in different order).
 Observe that in both octagons every edge goes from one of $\{a, b,c,d\}$ to one of $\{m,n,o,p\}$, but the vertices are shuffled in such a way that no edge is repeated.
 See Figure~\ref{fig:mani-walkup-pyr-o}.

\begin{figure}[hbt]
\begin{center}
\includegraphics[scale=.45]{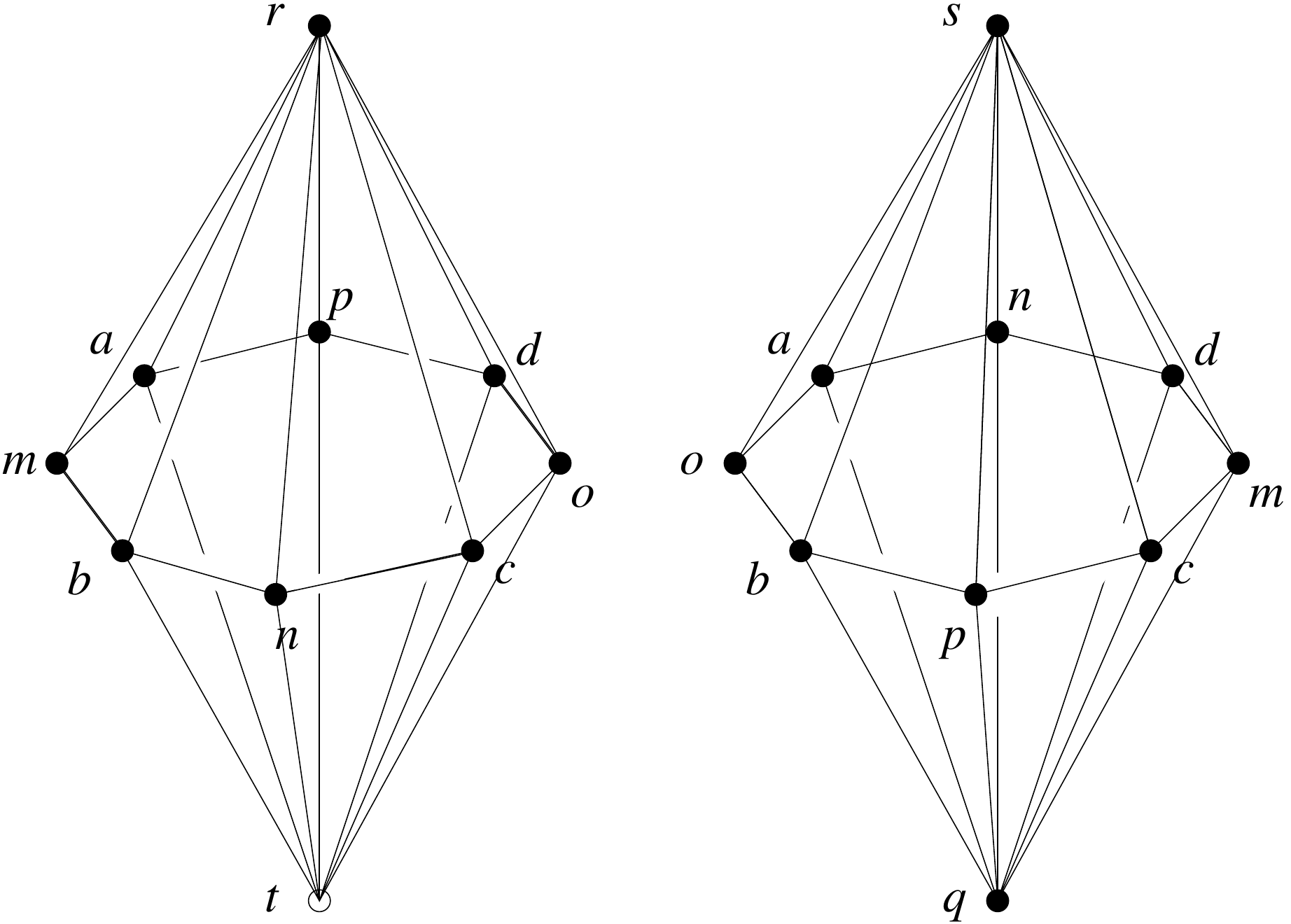}
\caption{Two octagonal bipyramids}
\label{fig:mani-walkup-pyr-o}
\end{center}
\end{figure}

The interiors of the two bipyramids can easily be triangulated (subdivided into terahedra) in such a way that the tetrahedron $abcd$ is used in the first one and $mnop$ in the second. Then the two bipyramids can be embedded in the $3$-sphere (with corresponding vertices identified) by first embedding them disjointly and then pinching the vertices of one of the octagons to glue them with those of the other. We claim that no extension of this partial triangulation to the whole $3$-sphere can have the non-revisiting property.

Indeed, every path from the tetrahedron $abcd$ to the tetrahedron $mnop$ must exit the first bipyramid through one of its boundary triangles, which uses one of the edges of the first octagon. In particular, our path will at this point have abandoned three of the vertices of $abcd$ and be using one of $mnop$. To keep the non-revisiting property, the abandoned vertices should not be used again, and the new one should not be abandoned, since it is a vertex of our target tetrahedron. But then it is impossible for our path to enter the second bipyramid: it should do so via a triangle using an edge of the second octagon, and non-revisiting implies that this edge should be the same used to exit the first bipyramid. This is impossible since the two octagons have no edge in common.

We skip the technical part of the proof, namely that $K$ can be completed to a triangulation of the $3$-sphere using 
the tetrahedra $abcd$ and $mnop$ (and with only four extra vertices). The way Mani and Walkup show it is by listing the tetrahedra of the whole triangulation and verifying that they form a shellable sphere. 
%The first claim, except for the number of extra vertices needed, follows from Whitehead's Completion Lemma (see~\cite{Whitehead}):\footnote{I am not sure this reference is correct. Perhaps a reference to Bing's book is better. paco} every simplicial complex embedded in a 3-manifold can be extended to a triangulation of the whole manifold.
%The number of extra vertices cannot \emph{a priori} be controlled, but that will only affect the number of vertices of the final 3-sphere  (and the number of wedges needed to get a non-Hirsch sphere of a certain dimension from it).
\end{proof}

%\subsection*{Acknowledgments}
%We thank David Bremner, Jes\'us De Loera, Antoine Deza, J\"org Rambau , G\"unter M. Ziegler  and the anonymous refereefor their valuable input. This paper grew out of several conversations during the second author's sabbatical leave at UC Davis in 2008 and the authors' participation in the I-Math DocCourse on Discrete and Computational Geometry at the \emph{Centre de Recerca Matem\`atica} in 2009. We thank both institutions for hosting us and the financial support from the National Science Foundation and  the Spanish Ministry of Science.  

\vskip.5cm
\noindent {\small Edward D. Kim}\newline
\emph{Department of Mathematics}\newline
\emph{University of California, Davis. Davis, CA 95616, USA}\newline
\emph{email: }\url{ekim@math.ucdavis.edu}\newline
\emph{web: }\url{http://www.math.ucdavis.edu/~ekim/}

\vskip.5cm
\noindent {\small Francisco Santos}\newline
\emph{Departamento de Matem\'aticas, Estad\'istica y Computaci\'on}\newline
\emph{Universidad de Cantabria, E-39005 Santander, Spain}\newline
\emph{email: }\url{francisco.santos@unican.es}\newline
\emph{web: }\url{http://personales.unican.es/santosf/}

\clearpage

\tableofcontents
\end{document}